\documentclass[10pt,letterpaper]{elsarticle}
\usepackage{graphicx}
\usepackage{epsfig}
\usepackage{amsmath}
\usepackage{amsfonts}
\usepackage{subfigure}
\usepackage{ulem}
\usepackage{bm}
\usepackage{color}
\newcommand{\vx}{\boldsymbol{x}}
\newcommand{\vp}{\boldsymbol{p}}
\newcommand{\vl}{\boldsymbol{p}}
\usepackage{ulem}

\newcommand{\xib}{{\bm \xi}}

\date{}

\begin{document}
\title{An Iterative Action Minimizing Method for Computing Optimal Paths in Stochastic Dynamical Systems}

\author{Brandon S. Lindley, Ira B.~Schwartz}

\address{U.S. Naval Research Laboratory\\ Code 6792, Plasma Physics Division, Nonlinear Systems Dynamics Section\\ Washington, D.C., 20375, USA\\email:brandon.lindley.ctr@nrl.navy.mil\\tel:202-404-8358 fax:202-767-0631}


\date{}
\begin{abstract}
We present a numerical method for computing optimal transition pathways and transition rates in systems of stochastic differential equations (SDEs). In particular, we compute the most probable transition path of stochastic equations by minimizing the effective action in a corresponding deterministic Hamiltonian system. The numerical method presented here involves using an iterative scheme for solving a two-point boundary value problem for the Hamiltonian system. We validate our method by applying it to both continuous stochastic systems, such as nonlinear oscillators governed by the Duffing equation, and finite discrete systems, such as epidemic problems, which are governed by a set of master equations. Furthermore, we demonstrate that this method is capable of dealing with stochastic systems of delay differential equations.
\end{abstract}
\maketitle




\section{Introduction}

One important aspect of the study of dynamical systems is the study of
  noise on the underlying deterministic
  dynamics~\cite{gar03,vanKampen_book}. Although one might expect
  the deterministic dynamics would be only slightly perturbed in the presence of small
  noise, there are now many examples where noise causes a dramatic measurable
change in behavior, such as noise induced switching between attractors in
continuous systems \cite{FW84} , and
noise induced extinction in finite size systems \cite{vanKampen_book}.

In systems transitioning between coexisting stable states, much research has been done primarily because switching can be now
investigated for a large variety of well-controlled micro- and mesoscopic
systems, such as trapped electrons and atoms,
Josephson junctions, and nano- and micro-mechanical oscillators
\cite{Lapidus1999,Siddiqi2005,Aldridge2005,Kim2005,Gommers2005,Stambaugh2006,Abdo2007,Lupascu2007,Katz2007,Serban2007}. In
these systems, observed fluctuations
are usually due to thermal or externally applied 
noise. However, as systems become smaller, an increasingly important role
may be played also by non-Gaussian noise. It may come, for example,
from one or a few two-state fluctuators hopping at random between the states, in which case
the noise may be often described as a telegraph noise. It may also be induced
by Poisson noise \cite{billings2010switching}.

In finite size populations or systems, extinction occurs in discrete, finite
  populations undergoing stochastic effects due to random transitions or
  perturbations.  The origins of stochasticity may be internal to the system  or may arise from the external environment~\cite{DeCastro2005,Lloyd2007}, and in most cases is non-Gaussian~\cite{schwartz2011converging,ForgostonBMB}.

Extinction depends on the nature and strength of the noise~\cite{Melbourne2008},  outbreak amplitude~\cite{Alonso2006} and seasonal phase
occurrence~\cite{stolhu07}.  For large populations, the intensity of
internal population noise is generally small.  However, a rare,
large fluctuation can occur with non-zero probability and the
system may be able to reach the extinct state.   Since the
extinct state is absorbing due to effective stochastic forces, eventual extinction is guaranteed when
there is no source of reintroduction~\cite{bartlett49,Allen2000,gar03}.

 Models of finite populations, which include extinction processes, are
effectively described using the master
equation formalism, and predict probabilities of rare events~\cite{schwartz2011converging}. For many problems involving extinction in large populations, if the
probability distribution of the population is quasi-stationary, the probability of extinction
is a function that decreases exponentially with increasing population size.  The exponent in this function scales as a
deterministic quantity called the action~\cite{Kubo1973}.  It can be shown that a
trajectory that brings the system to extinction is
very likely to lie
along a most probable path, called the  optimal path.  It is a major property that a
deterministic quantity such as the action can predict the probability of
extinction, which is inherently a stochastic process, and is also formulated
in continuous systems driven by noise~\cite{sbdl09,dyscla08}.

Locating the optimal path is important since the quantity of interest,
whether it is the switching or extinction rate, depends on the probability to
traverse this path.  Therefore,  a stochastic control strategy based on the switching or
extinction rates can be determined by its effect on the optimal
path~\cite{dyscla08}.

 The  optimal path formalism converts the entire
stochastic problem to a mechanistic dynamical 
systems problem with definitive properties. First, the optimal path is a
solution to a Hamiltonian dynamical system. In the case of continuous
stochastic models, the dimension of the system is
twice that of the original stochastic problem. The other dimensions are
conjugate momenta, and typically represent the physical force of the noise
which induces escape from a basin of attraction to either switch or go
extinct. Finally, due to the symplectic structure of the resulting Hamiltonian
system, it can be shown that both attractors and saddles of the original system become saddles of the Hamiltonian system.

One of the main obstacles to finding the optimal path is that it is an
  inherently unstable object. That is, if one one starts near the path described
  by the Hamiltonian system, then after a short time, the dynamics leaves the
  neighborhood of the path. In addition, although the path may be hyperbolic near
  the saddle points, it may not be hyperbolic along the rest of the
  path. Solving such problems using shooting methods for simple epidemic
  models~\cite{keller1992numerical,Meerson2008,Meerson2012}, or mixed shooting using forward and
  backward iteration~\cite{Elgart2004,Chernykh2001} will in general be inadequate to handle
  even the simplest unstable paths in higher dimensions.
  Therefore, it is the goal of this paper to exemplify a robust numerical
method to solve for the optimal path using a general accurate discrete
formulation applied to the Hamiltonian two point  boundary value problem.

The method we employ here to compute the optimal paths is similar to the generalized minimum action method (gMAM), \cite{Heymann2008} which is a blend of the string method \cite{WeinanE2002} and minimum action method \cite{WeinanE2004}. Both are iterative methods which globally minimize the action along the path. These other techniques differ from ours primarily in that our formulation of the problem allows a direct, fully explicit iterative scheme, while the gMAM, in particular, employs a semi-implicit scheme. The numerical scheme presented here should provide an easy to employ alternative to the methods discussed above for stochastic optimization problems which are formulated as Hamiltonian two-point boundary value problems.

The paper is organized as follows. We first briefly present the general SDE problem, and the formulation of the corresponding deterministic Hamiltonian system by treating the switching as a rare event. We then present the details of the numerical approximation technique we will use to find the path which maximize the probability of switching. Using this technique, we then demonstrate finding the optimal path from a stable focus to the unstable saddle for the unforced Duffing equation, then compute optimal extinction pathways for a simple epidemic model, and finally adapt our method to find optimal transition paths in stochastic delay differential systems.

\section{General Problem}

Consider a general stochastic differential equation of the form
\begin{equation}
\dot{\vx}(t) = \boldsymbol{f}(\vx(t)) + \boldsymbol{G}(\vx(t))\boldsymbol{\xi}(t),
\label{general:ode}
\end{equation}
where $\boldsymbol{x}\in\mathbb{R}^n$ represents the physical quantity in
  state space and
 and the matrix $\boldsymbol G$\footnote{Throughout the paper, boldface lower-case letters will indicate vectors, while boldface upper-case letters will indicate matrices. } is given by
$\boldsymbol{G}(\vx(t))=\text{diag}\{g_1(t),g_2(t),...g_n(t)\}$, where the $g_i$'s are general nonlinear functions. We suppose the
  noise  $\boldsymbol{\xi}\in\mathbb{R}^n$ is a  vector having a
  Gaussian distribution with intensity $D$, and independent components. It is characterized by its probability density
  functional 
  $\mathcal{P}_{\boldsymbol{\xi}}=e^{-\mathcal{R}_{\xi}/D}$,  
\begin{equation}
\label{eq:Gauss_Functional}
{\cal R}_{\xib}[\xib(t)]= \frac{1}{2}\int dt\,dt'\,\xib(t)\xib(t').
\end{equation}
 We wish to determine the path with the maximum probability of traveling
from the initial state $\bm{x}_A$ to the final state $\bm{x}_B$, where the
initial and final states are equilibria of the noise-free (i.e. $\boldsymbol \xi = \boldsymbol 0$) version of
Eq. \ref{general:ode} given by $\bm{f(x)}=\bm{0}$. These states typically
  characterize a generic problem of study in stochastic systems, such as switching between
  attractors, escape from a basin of attraction, or extinction of population. We assume the noise intensity $D$ is sufficiently small so
that in our analysis sample paths will limit on an optimal path as $D
\rightarrow 0$. We also remark that $\bm{\xi}$ is formally the time derivative
of a Brownian motion, sometimes referred to as white noise \cite{fleming1975deterministic}.

For $D$ sufficiently small, and examining the tail of the distribution for a large
fluctuation (which is assumed to be a rare event), the probability of observing such a
large fluctuation  scales exponentially by \cite{Freidlin_book,Dykman1990},

\begin{equation}
\begin{aligned}
\mathcal{P}_{x} = e^{-R/D}, && R = \text{min}\mathcal{R}(\vx,\boldsymbol{\xi},\boldsymbol{p}),
\end{aligned}
\label{rare:event}
\end{equation}
where,
\begin{equation}
\mathcal{R}(\vx,\boldsymbol{\xi},\boldsymbol{p}) = R_{\xi}+\int
\boldsymbol{p}\cdot [
\dot{\vx}-\boldsymbol{f}(\vx)-\boldsymbol{G}(\vx(t))\boldsymbol{\xi}
], 
\label{exponent}
\end{equation}
where the Lagrange multipliers $\boldsymbol{p}$ also correspond to the
conjugate momenta of the equivalent Hamilton-Jacobi formulation of this problem.\footnote{The vector multiplication here is assumed to be an inner product.} The exponent
$R$ of Eq. \ref{rare:event} is called the action, and corresponds to the minimizer of the action in the Hamilton-Jacobi formulation which occurs along the optimal path. This path
will minimize the integral of Eq.~\ref{exponent}, and is found by setting
the variations along the path $\delta\mathcal{R}$ to zero.

 The resulting equations of motion for the states and Lagrange multipliers
  are given by 
\begin{equation}
\label{eq:EOM}
\begin{aligned}
\dot{\vx} &= \boldsymbol{f}(\vx(t)) + \boldsymbol{G}^2(\bm{x})\vl &&\\
\dot{\vl} &= -\bm{G(x)}{\frac{\partial \bm{G}}{\partial \bm{x}}(\bm{x})}\vl\vl
-\frac{\partial \bm{f(x)}}{\partial \bm{x}}\vl.
\end{aligned}
\end{equation}
Here $[\frac{\partial \bm{G}}{\partial
    \bm{x}}(\bm{x})\vl\vl]_i=[\frac{\partial \bm{G}}{\partial
    \bm{x}}(\bm{x})]_{ijk}[\vl]_j[\vl]_k$, where Einstein summation is assumed
over repeated indices. Note that $\vl=0$ is invariant, and recovers the noise free case of Eq. \ref{general:ode}. The above equations can be shown to satisfy the motion of a Hamiltonian system
with Hamiltonian
\begin{equation}
H(\bm{x},\vl) =
\frac{(\boldsymbol{G}^2(\boldsymbol{x})\boldsymbol{p})\cdot\boldsymbol{p}}{2}+\boldsymbol{p}\cdot
\boldsymbol{f}(\vx).\label{hamiltonian}
\end{equation}
That is, the dynamics of a given path satisfy $\dot{\vx} = \frac{\partial
  H(\bm{x},\bm{p})}{\partial \vp}$ and  
$\dot{\vp} = -\frac{\partial H(\bm{x},\bm{p})}{\partial \vx}.$


In addition to solving the Hamiltonian system of dynamics, the full
  problem specification of an optimal path requires  boundary conditions for
  both state $\bm{x}$ and momenta $\vp$. The boundary conditions of the
  optimal path consist of two steady states,
$X_A=(\boldsymbol{x}_A,\boldsymbol{p}_A)$ and
$X_B=(\boldsymbol{x}_B,\boldsymbol{p}_B)$ at equilibrium. Typically, the boundary conditions
are derived from the equations of motion defined  by  Eq.~\ref{hamiltonian}. Since they
are derived as steady state conditions, they are asymptotic boundary
conditions that are infinite limits in the temporal line. In addition, since
we have assumed the Hamiltonian is time-independent, energy is conserved, and
the path must lie on a fixed energy surface. Notice, in the case where the Hamiltonian is time invariant, and given that the action is minimized along the optimal path, the Hamilton-Jacobi equations require a zero-energy constraint, i.e. $H(\vx,\vp)=0$. We now describe how to solve the
Hamiltonian system as a two point boundary value problem on a restricted
energy surface.

\subsection{Stability of Steady State Solutions}

As stated above, we seek the optimal path between two steady state solutions
$X_A=(\boldsymbol{x}_A,\boldsymbol{p}_A)$ and
$X_B=(\boldsymbol{x}_B,\boldsymbol{p}_B)$ where $\boldsymbol{x}_A$ and
$\boldsymbol{x}_B$ are steady state solutions of the zero-noise case of
Eq. \ref{general:ode}, and thus satisfy $\boldsymbol{f}(\boldsymbol{x}_A)=
\boldsymbol{f}(\boldsymbol{x}_B)=0$. Depending upon the linear stability of
the zero-noise steady states, the stochastic transition from $\vx_A$ to
$\vx_B$ may or may not be a rare event. For example, if the stability matrix
of the general SDE Eq. \ref{general:ode} $\boldsymbol{A}=\frac{\partial
  \boldsymbol{f}(\vx_A)}{\partial\vx}=\boldsymbol{f}'(\vx_A)$ has at least one
eigenvalue with positive real part, then it is an unstable or saddle point in
the stochastic equation. Further, if  $\boldsymbol{f}'(\vx_B)$ has eigenvalues
with all negative real parts, then it is a stable focus, and the transition
from $\vx_A$ to $\vx_B$ will be deterministic restricted to $\vl = \bm{0}$ and therefore not a rare event. The rare event, in this case, would be a transition from $\vx_B$ to $\vx_A$, meaning noise drives the system out of the stable focus and onto the unstable/saddle point.

It is worth noting that all steady state solutions of the Hamiltonian system \ref{eq:EOM} which correspond to the deterministic steady state solutions of Eq. \ref{general:ode}, i.e. $\boldsymbol{f}(\vx)=0$, are saddle points. It is easy to see why this is true in the case of additive noise, i.e. $\boldsymbol{G} = \boldsymbol{I}$. We can classify the steady states of Eq. \ref{eq:EOM} by calculating the eigenvalues of its stability matrix $\boldsymbol{Q}$
\begin{equation}
\boldsymbol{Q}=\begin{bmatrix}
\boldsymbol{f}'(\vx_{A,B}) && \boldsymbol{I} \\ -\boldsymbol{f}''(\vx_{A,B})\vp_{A,B}  && -\boldsymbol{f}'(\vx_{A,B}) 
\end{bmatrix}=\begin{bmatrix}
\boldsymbol{A} && \boldsymbol{I} \\ \boldsymbol{0}  && -\boldsymbol{A} ,
\end{bmatrix},
\end{equation}
at the steady states. The eigenvalues $\gamma_k$ of $\boldsymbol{Q}$ are given by,
\begin{equation}
\begin{bmatrix}
\boldsymbol{A} && \boldsymbol{I} \\ \boldsymbol{0}  && -\boldsymbol{A} ,
\end{bmatrix} \begin{bmatrix}\boldsymbol{z}_1\\ \boldsymbol{z}_2\end{bmatrix}=\begin{bmatrix}
\boldsymbol{A}\boldsymbol{z}_1 + \boldsymbol{z}_2 \\ -\boldsymbol{A}\boldsymbol{z}_2
\end{bmatrix}=\gamma_k  \begin{bmatrix}\boldsymbol{z}_1\\ \boldsymbol{z}_2\end{bmatrix}.
\end{equation}
If we assume $\boldsymbol{z}_2=0$, then
$\boldsymbol{A}\boldsymbol{z}_1=\gamma_k\boldsymbol{z}_1$. Thus, if
$\lambda_k\in\mathbb{R}$ is an eigenvalue of $\boldsymbol{A}$, the stability
matrix of the SDE problem given by Eq. \ref{general:ode}, then
$\gamma_k=\lambda_k$ is an eigenvalue of $\boldsymbol{Q}$ with eigenvector
$[\boldsymbol{z}_1,\boldsymbol{0}]^T$. Similarly, if  $\boldsymbol{z}_1=0$,
then -$\boldsymbol{A}\boldsymbol{z}_2=\gamma_k\boldsymbol{z}_2$, which implies
that  $\gamma_k=-\lambda_k$ is an eigenvalue of $\boldsymbol{Q}$ with
eigenvector $[\boldsymbol{0},\boldsymbol{z}_2]^T$. Since $\pm\lambda_k$ are
eigenvalues of $\boldsymbol{Q}$, we can conclude that every steady state
solution of Eq. \ref{eq:EOM} has eigenvalues with both positive and negative
real parts and, thus, is a saddle. Thus, every steady state solution whose
linearization has non-zero real part of the general stochastic equation,
Eq. \ref{general:ode} regardless of stability becomes a saddle point when the
system is converted to a Hamiltonian system. That is, deterministic
  attractors and saddles map to saddles in the Hamiltonian formulation.

\section{Numerical Scheme}
Our numerical approach involves using a finite differences scheme to write the system of ODEs as a high dimensional algebraic system, to which we apply a modified Newton's Method to minimize the residual error until a solution is reached (to within some desired tolerance). As before, assume the Hamiltonian system in $\mathbb{R}^{2n}$ admits two steady states, $X_A$ and $X_B$. 

Then we seek the optimal path on the zero-energy surface that connects $X_A$ to $X_B$. In this formulation, one would expect such a path to exhibit several properties when parametrized along $t\in(-\infty,\infty)$. First, we assume the path starts at $X_A$ at $t=-\infty$. Since this point is an equilibrium solution, we should expect that the solution stays very near this value, that is, there exists $\epsilon>0$ such that $|X_A-X(t)|<\epsilon$ for $-\infty< t\leq -T_\epsilon$, has a transition region from $-T_\epsilon< t < T_\epsilon$ and finally stays near $X_B$, the second steady state for, $T_\epsilon\leq t < \infty$. Numerically, we will approximate the solution on the finite domain $-T_\epsilon\leq t \leq T_\epsilon$ so that the value at $X(\pm T_\epsilon$) is arbitrarily close to the steady solution. 

We map the interval $[-T_\epsilon,T_\epsilon]$\footnote{For the simulations below, $|T_{\epsilon}|\geq 100$ unless otherwise specified. Ideally, $T_{\epsilon}$ is picked large enough so that $|\vx_A-\vx(-T_\epsilon)|<10^{-16}$ and $|\vx_B-\vx(T_\epsilon)|<10^{-16}$. i.e. the steady states are obtained up to machine precision for double precision numbers.} onto $[0,1]$ using the linear transformation $t = 2T_{\epsilon}\bar{t}-T_{\epsilon}$, and drop the ``bar'' notation for readability. On this discrete time domain, we write the system of ODEs as a system of nonlinear algebraic equations using central differences. 

The simplest method is to employ a finite step size $h = 1/N$ and use a uniform time step subdividing $[0,1]$ into $N+1$ equal segments. In practice, however, the simple uniform step size is not always the best choice since the optimal path tends to stay very near the stable points throughout most of the domain, and sometimes makes a relatively sharp transition near the center of the domain. In this case, it is helpful to use a nonuniform grid to resolve the sharp transition region using a fine mesh, and to use a coarse mesh near the edges where the solution is mostly flat. Thus, for the nonuniform time step $h_k$, yielding the time series $t_{k+1}=t_k+h_k$ and corresponding function values $\vx_k$, the derivative is approximated by the operator $\delta_{h}$,

\begin{equation}
\frac{d}{dt}\vx_k \approx \delta_h \vx_k\equiv \frac{h_{k-1}^2\vx_{k+1} + (h_{k}^2-h_{k-1}^2)\vx_k-h_{k}^2\vx_{k-1}}{h_{k-1}h_{k}^2+h_{k}h_{k-1}^2}.
\end{equation} 

At this point, we can write the generic system of $2n(N+1)$ nonlinear algebraic equations:

\begin{equation}
\begin{aligned}
\delta_h\vx_k- \frac{\partial
  H(\vx_k,\vp_k)}{\partial \vp}=0 && \delta_h \vp_k+ \frac{\partial
  H(\vx_k,\vp_k)}{\partial \vx}=0, && k = 0,1,...,N,
\label{discrete:equations}
\end{aligned}
\end{equation}
and solve this system using a general Newton's Method for nonlinear systems of
equations. To properly apply Newton's method here, let $\bm{q}_j =
  \{\vx_{1,j}...\vx_{N,j},\vp_{1,j}...\vp_{N,j}\}^T$ be the extended vector of
  dimension $1\times2nN$ containing the j-th Newton iterate (recalling that
  $\vx_{k,j},\vl_{k,j}\in\mathbb{R}^n$ are defined on the timeseries given by
  $k=0,1,...,N$). Then $j=0$ will represent the initial guess. Let
  $\mathcal{F}:\mathcal{R}^{2nN}\rightarrow\mathcal{R}^{2nN}$ be the
  function defined by Eq. \ref{discrete:equations} acting on $\bm{q}_j$. Then
to find the zeros of $\mathcal{F}(\bm{q})$ we employ a Newton scheme. A new Newton iterate is given by solving the linear system $J_{\mathcal{F}}(\bm{q}_j)(\bm{q}_{j+1}-\bm{q}_j) =- \mathcal{F}(\bm{q}_j)$, using any one of a variety of methods such as LU decomposition or the generalized minimal residual method (GMRES) with appropriate preconditioners. Throughout this paper we will use LU decomposition with partial pivots optimized for a sparse linear system. Here the Jacobian $J_{\mathcal{F}}(\bm{q}_n)$ is computed approximately using a central difference scheme.

Formally, this method is second order with respect to $h_k$. The initial guess for this algorithm is constructed by the knowledge that the optimal path spends most of its time near the stable equilibria, and has a brief but sometimes sharp transition between the two states. One choice that has worked in practice is using functions like $\bm{q}_{0,k}=(\vx_A - \vx_B)/(1+e^{Ct_k})+\vx_B$, with $k=0,1,...N$ (where the $C>0$ parameter adjusts the sharpness of the jump), which have horizontal asymptotes at the appropriate critical values. Usually, though not always, $\bm{q}_{0,k}$ is set so that for $k = N+1,N+2,...2N$, $\bm{q}_{0,k}=\bm{0}$.

Note that the zero-energy surface constraint $H(\boldsymbol{x},\boldsymbol{p}) = 0$ is not imposed. Rather, the initial guess will start out near $X_A$ (at $t=0$) and $X_B$ (at $t=1$) which lie asymptotically close to the zero energy surface, and thus the final solution will have to lie on this surface since such a solution is time invariant (i.e. $d/dt H(\vx,\vp)=0$). At each iterate, both the residual error and the Hamiltonian are checked at each point, and both must reach a desired tolerance in the $\mathcal{L}_\infty$ norm before the procedure is completed.

Once the optimal path is computed, the action (i.e. the exponent) along the optimal path may be obtained with a simple integral,
\begin{equation}
R = \int_{t=0}^1 \left[\vp(t)\cdot\frac{d\vx(t)}{dt}-H(\vx(t),\vp(t))\right]dt.
\label{action:integral}
\end{equation}
Exhaustive convergence tests on the residual error for a variety of test problems for both the uniform and non-uniform grids have demonstrated the second order convergence for this method. We have noted some dependency on the initial guess in terms of the overall speed of convergence (or divergence for a particularly bad guess). The method generally produces a unique solution (up to possibly a horizontal shift in the time series seen in a few examples, which does not affect the path integrals of interest). This method has been reliable for a wide parameter regime for each of the test problems, with the limitations to be discussed below on a case by case basis. 

This method, which we will henceforth refer to as the Iterative Action Minimization Method (IAMM), has several distinct advantages over other methods, the foremost of which is straightforward scalability to higher dimensions. For very high dimensional problems, the systems will eventually become too large to treat easily with a single processor, but this algorithm has proven efficient for up to six dimensional problems with a single processor. Further, this method lends itself to infinite dimensional problems, such as time delay stochastic differential equations, as will be demonstrated below.

\section{Noise Induced Transitions}
We demonstrate the numerical techniques by examining several bistable dynamical systems. Using methods discussed above, we explicitly approximate the optimal path between the two states and then numerically integrate along the path directly to compute the action. 
\subsection{Switching in the Duffing Equation}
One of the standard nonlinear dynamical systems which exhibits bi-stability is Duffing's equation. This equation is used to model certain types of nonlinear damped oscillators, and here we consider the singularly perturbed and unforced version \cite{forsch09},
\begin{equation}
\begin{aligned}
\dot{x} &= y + g_1(x,y)\xi_1(t),\\ 
\epsilon\dot{y} &= \alpha x - \beta x^3 - \delta y + g_2(x,y)\xi_2(t)\\
\dot{\epsilon} &= 0,
\end{aligned}
\label{2dDuff}
\end{equation}
Here $\alpha$ and $\beta$ control the size and nonlinear response of the restoring force, while $\delta$ controls the friction or damping on the system. The terms  $\xi_1$ and $\xi_2$ are uncorrelated white noise sources applied to the acceleration and velocity respectively. When the perturbation $\epsilon\ll 1$ fast and slow manifolds can be identified, while $\epsilon=1$ gives the unconstrained case. Rescaling time $t'=(1/\epsilon) t$, applying eq. \ref{hamiltonian}, and following the methodology above, we can write the system in the following general form,
\begin{equation}
H(x,y,p_x,p_y) = \frac{\epsilon p_x^2g_1(x)^2}{2} +  \frac{p_y^2g_2(y)^2}{2} + \epsilon p_x y + p_y (\alpha x - \beta x^3 - \delta y) = 0,
\label{DuffHam}
\end{equation}
with corresponding equations of motion,
\begin{equation}
\begin{aligned}
\dot{x} &= \epsilon(p_xg_1(x)^2 + y)\\ 
\dot{y} &= p_yg_2(y)^2 + \alpha x - \beta x^3 - \delta y\\ 
\dot{p}_x&=-\epsilon p_x^2g_1(x)g_1'(x) - p_y ( \alpha-3\beta x^2) \\  
\dot{p}_y&=-p_y^2g_2(y)g_2'(y) - \epsilon p_x + \delta p_y.
\end{aligned}
\label{4dDuff}
\end{equation}
 We will restrict our focus to the additive white noise case, $g_1(x,y)=g_2(x,y) = 1$.

Note that the Hamiltonian system emits three known steady states $X_A =
(0,0,0,0)$ and $X_B = (\pm\sqrt{\alpha/\beta},0,0,0)$, all of which are saddle
points. These steady states correspond to zero-noise critical points of
Eq. \ref{2dDuff}, $\vx_A=(0,0)$, a saddle point, and
$\vx_B=(\pm\sqrt{\alpha/\beta},0)$, the centers of the stable foci. The path
from $x_A$ to $x_B$ can be found deterministically when $\vl_x=\vl_y=0$, as
any solution perturbed from the saddle node point will move along the solution
curves and end up at either stable focus. A more interesting case is the
optimal path from one of the focus points to the saddle-node point, which will
require non-trivial momentum. Such momenta model the small noise effects which
organize to force the trajectory across the basin of attraction, thus escaping
from one attractor to the other. 

\begin{figure}[h!]
\begin{center}
\subfigure[]{\includegraphics[scale=0.35]{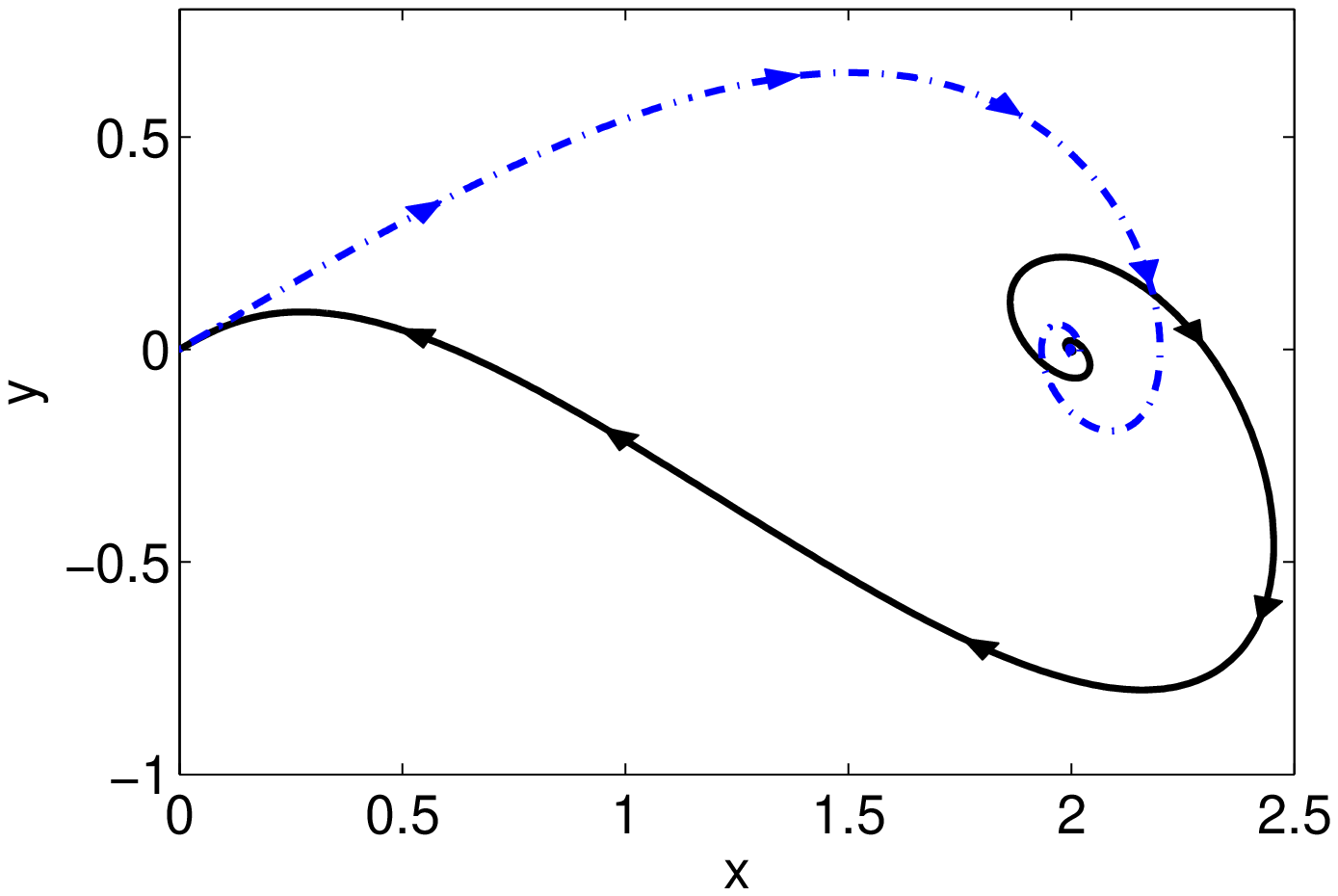}}
\subfigure[]{\includegraphics[scale=0.35]{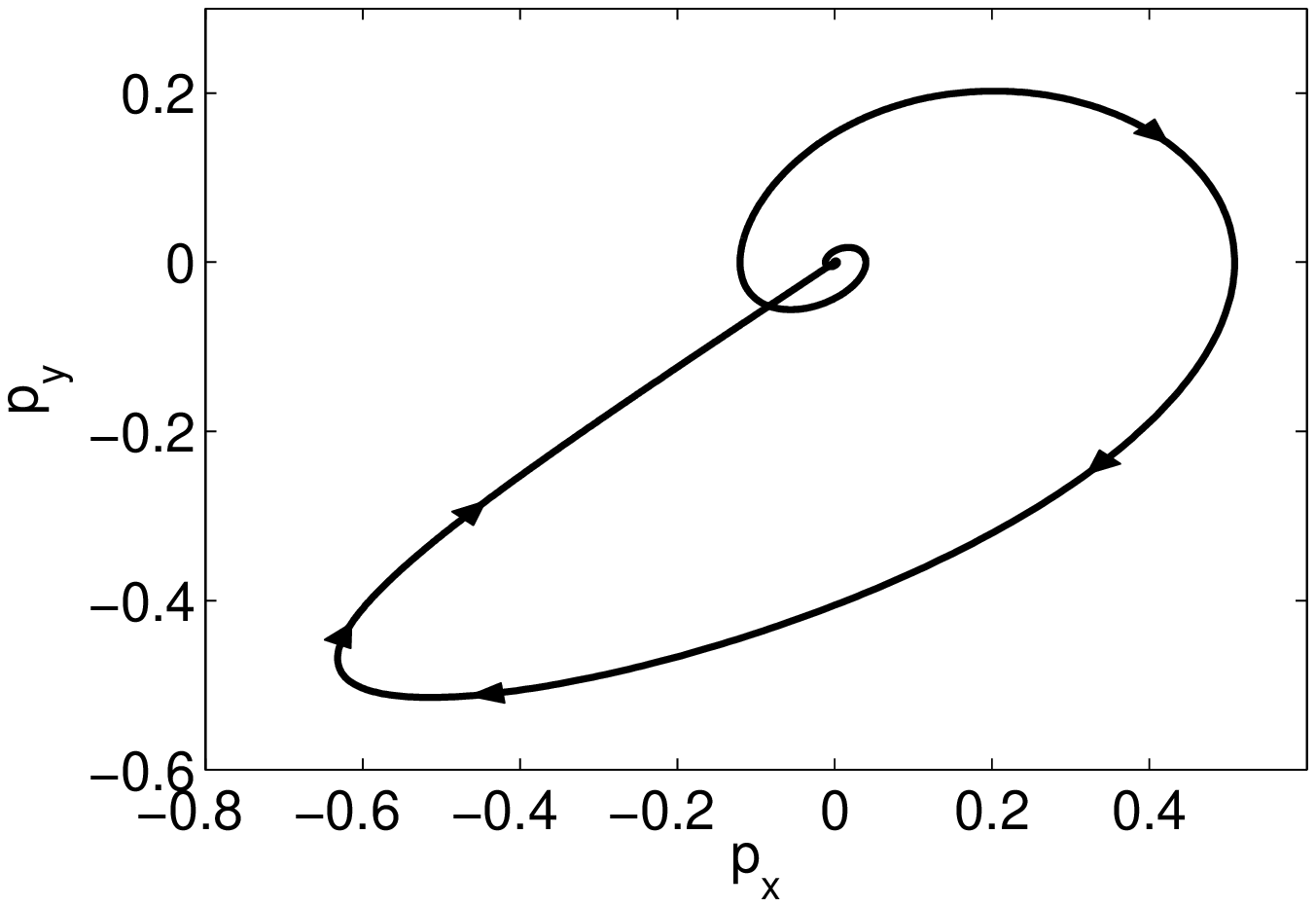}} 

\caption{Projections of the optimal path (a) onto phase space and (b) onto momenta space. In (a) we show both the optimal (from $(\sqrt{\alpha/\beta},0)$ to $(0,0)$, labeled with a solid line) and deterministic path  (from $(0,0)$ to $(\sqrt{\alpha/\beta},0)$, labeled with a dashed line) as predicted by our numerical method. Along the deterministic path, the action is zero and $p_x$ and $p_y$ are zero, but the corresponding noise components along the optimal path are shown in (b). Here $\epsilon=1$, $\beta=.25$, $\delta=1$, and $\alpha=1$.}\label{duff_bapts}
\end{center}
\end{figure}

\begin{figure}[h!]
\begin{center}
\subfigure[]{\includegraphics[scale=0.35]{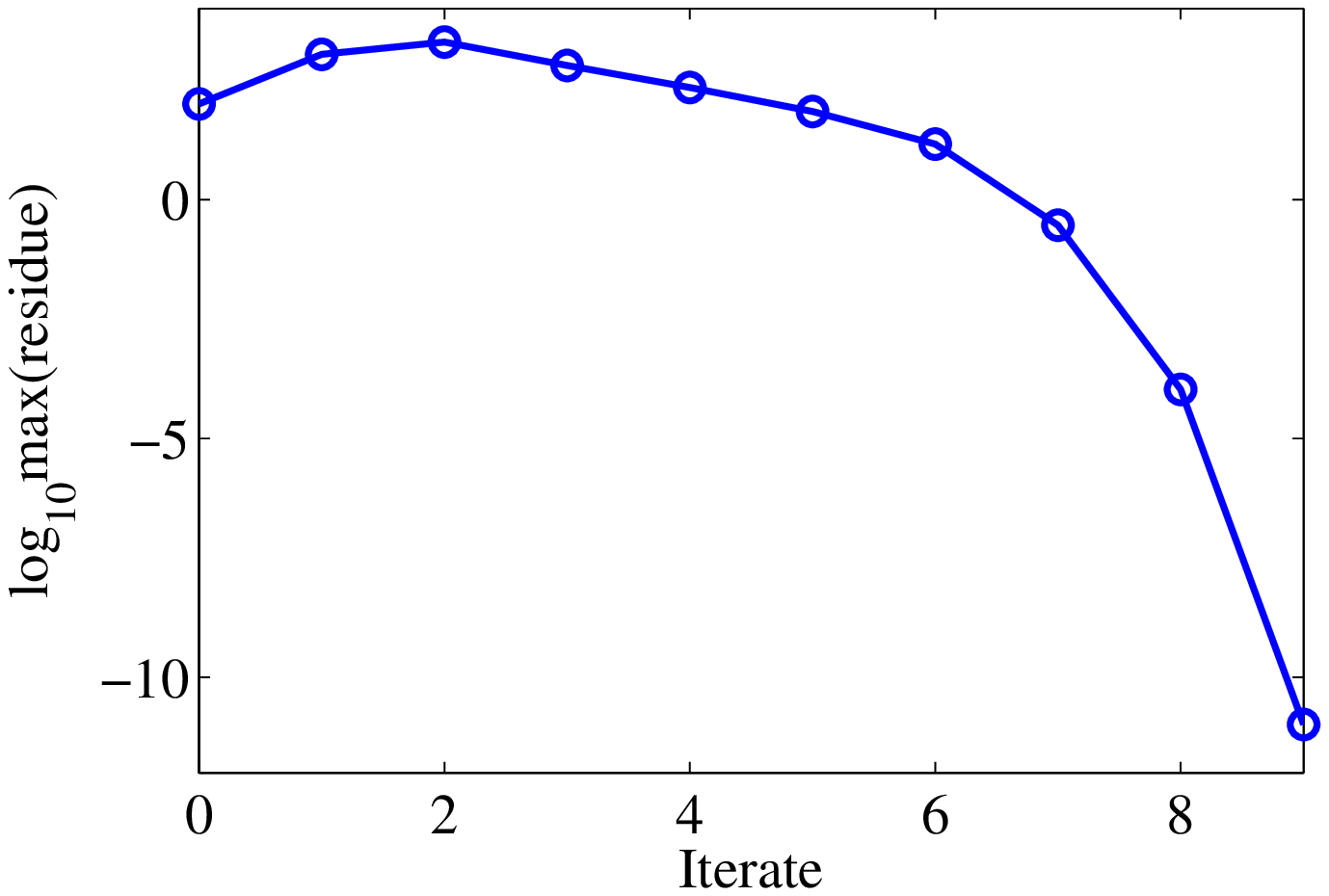}}
\subfigure[]{\includegraphics[scale=0.35]{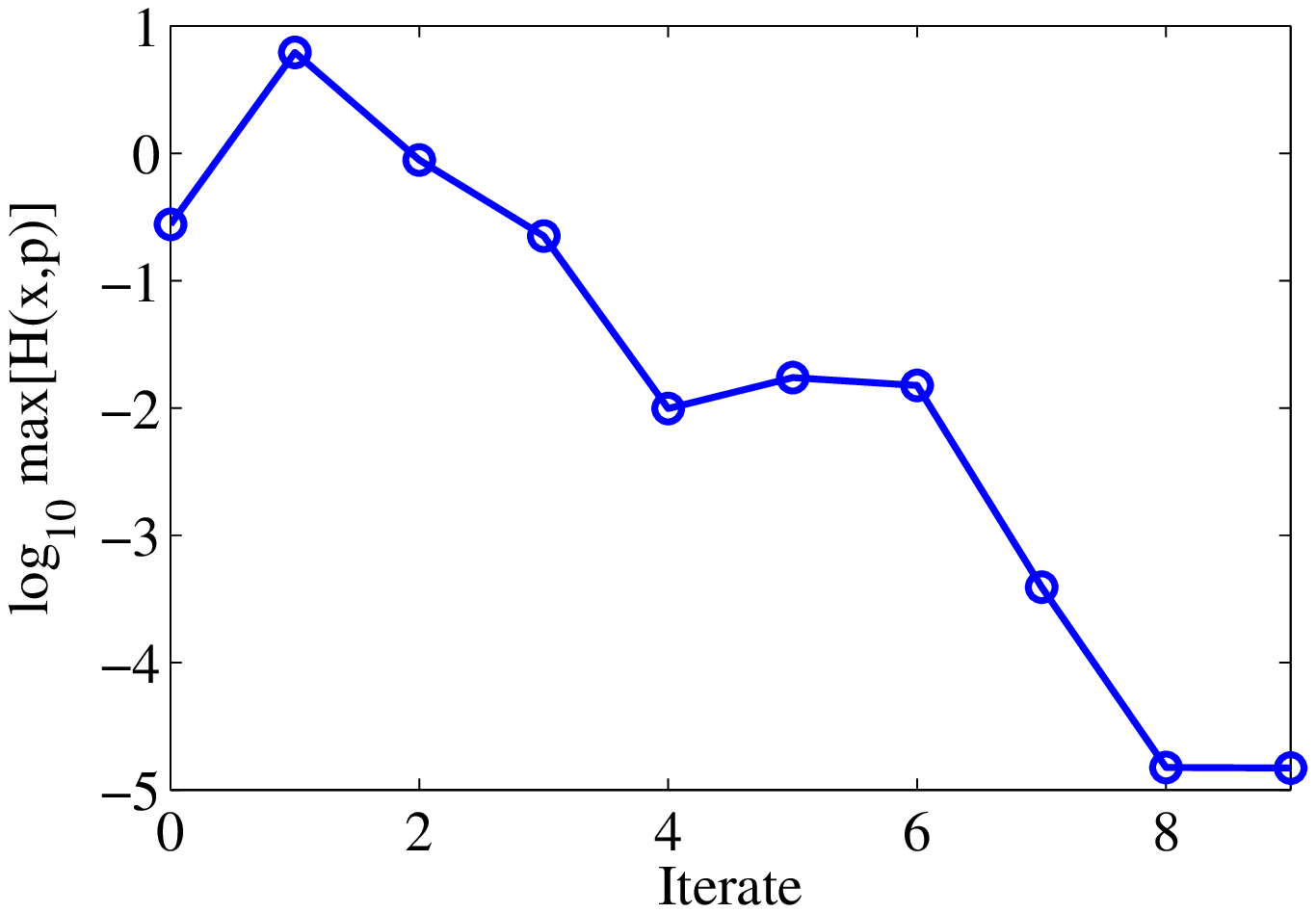}}
\caption{The maximal residual error (a) and the Hamiltonian (b) of the numerical scheme as a function of iterate, where the initial guess is iterate $0$. This figure shows the convergence for the data shown in Figure \ref{duff_bapts} where the parameter values are $\epsilon=1$, $\beta=.25$, $\delta=1$, and $\alpha=1$, and the error threshold is met after the ninth iteration. The maximal error is defined to be the maximum residue of all the time series for all the components.}\label{duff_converge}
\end{center}
\end{figure}

Figure \ref{duff_bapts} shows both the deterministic and optimal paths as computed using the IAMM
developed above. For the deterministic path, the algorithm correctly predicts
that the noise will be zero along this path, and that the action will be zero
as the probability of going from $X_A$ to $X_B$ is one, and thus not a rare
event. On the other hand, the path from $X_B$ to $X_A$ involves nontrivial
action, and the effect of the noise along the optimal path is
shown. This path lies on the zero-energy surface, and maximizes the
probability of traveling from $X_B$ to $X_A$ for arbitrarily small noise
  intensities, $D$. Figure \ref{duff_converge} shows the residual error at each iterate used to generate the data for figure \ref{duff_bapts} until the convergence criteria is reached. 

When $\epsilon\ll 1$, we can derive an analytical formulation for the action by using a center manifold analysis on Eq. \ref{DuffHam} by following the work of \cite{forsch09}.
Here the center manifold is given by $y=h(x,\epsilon)$, and we approximate it as,
\begin{equation}
y=h(x,\epsilon) = h_0(x) + \epsilon h_1(x) + \mathcal{O}(\epsilon^2).\label{centermani}
\end{equation}
By substituting equation \ref{centermani} into equation \ref{2dDuff} and equating like powers of $\epsilon$, we arrive at a one dimensional form of equation \ref{2dDuff} for the lowest order terms involving $\epsilon$,
\begin{equation}
\dot{x} = \frac{\epsilon}{\delta}\left(\alpha x-\beta x^3+\xi(t)\right).\label{1dDuff}
\end{equation}

Here the contribution of the uncorrelated noise terms $\xi_{1,2}$ are contained in a single noise source $\xi$. The Hamiltonian form of \ref{1dDuff} is,
\begin{equation}
H(x,p) = \frac{p^2}{2} + p\frac{\epsilon}{\delta}\left(\alpha x-\beta x^3\right).\label{1dDuff_Ham}
\end{equation}
We can find the nontrivial relationship ($p\neq 0$)between $p$ and $x$ on the zero energy surface directly, $p=-\frac{2\epsilon}{\delta}(\alpha x - \beta x^3)$ , and integrate along this path to predict the action along the optimal path. Since one may be interested in how the action scales relative to the distance between the two critical points we substitute the values $a^2=\alpha$ and $b^2=\frac{1}{\beta}$ to better illustrate the scaling with respect to $\beta$. Thus, with this substitution, we seek the action from $x=ab$ to $x=0$. From equation \ref{action:integral}, this is a simple integral,
\begin{equation}
R = -2\frac{\epsilon}{\delta}\int_{ab}^0 \left(a^2x-\frac{x^3}{b^2}\right) dx = \frac{\epsilon}{2\delta}a^4b^2.\label{action_scaling}
\end{equation}
Thus, we can predict, for example, that the action from $(ab,0)$ to $(0,0)$ should scale like the square of $b$ near the center manifold. Indeed, varying the parameter $b$ for the four-dimensional system, Eq. \ref{4dDuff}, predicts the same order scaling as seen below in Fig.~\ref{duff_center_scale}. We also consider the scaling with respect to the damping parameter $\delta$, and again, near the center manifold, the predictions are born out by integrating along the optimal path. Indeed, the scaling for the action predicted from the lowest order terms in the center manifold analysis seems to persist in the two-dimensional model.

\begin{figure}[h!]
\begin{center}
\subfigure[]{\includegraphics[scale=0.35]{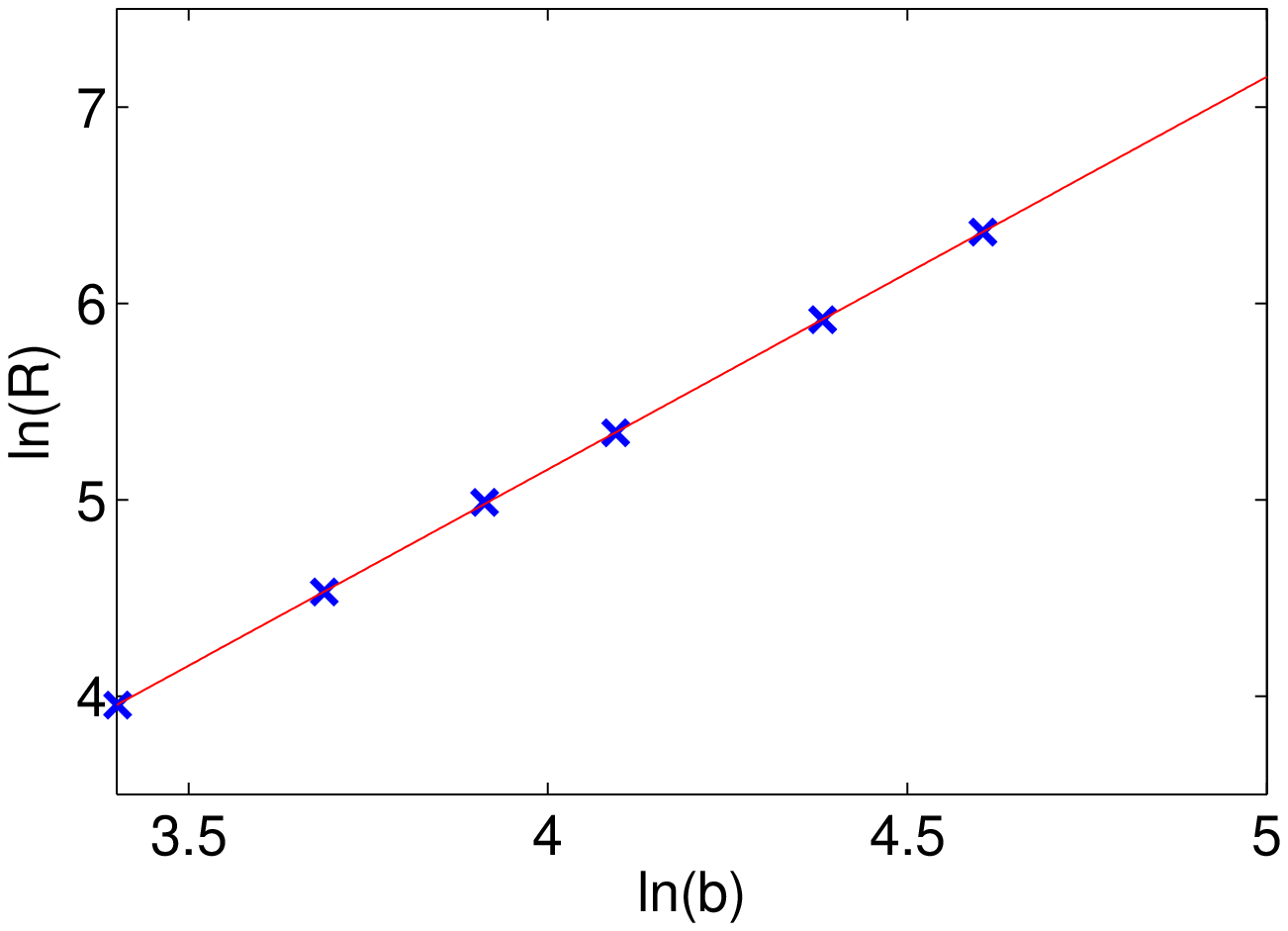}}
\subfigure[]{\includegraphics[scale=0.35]{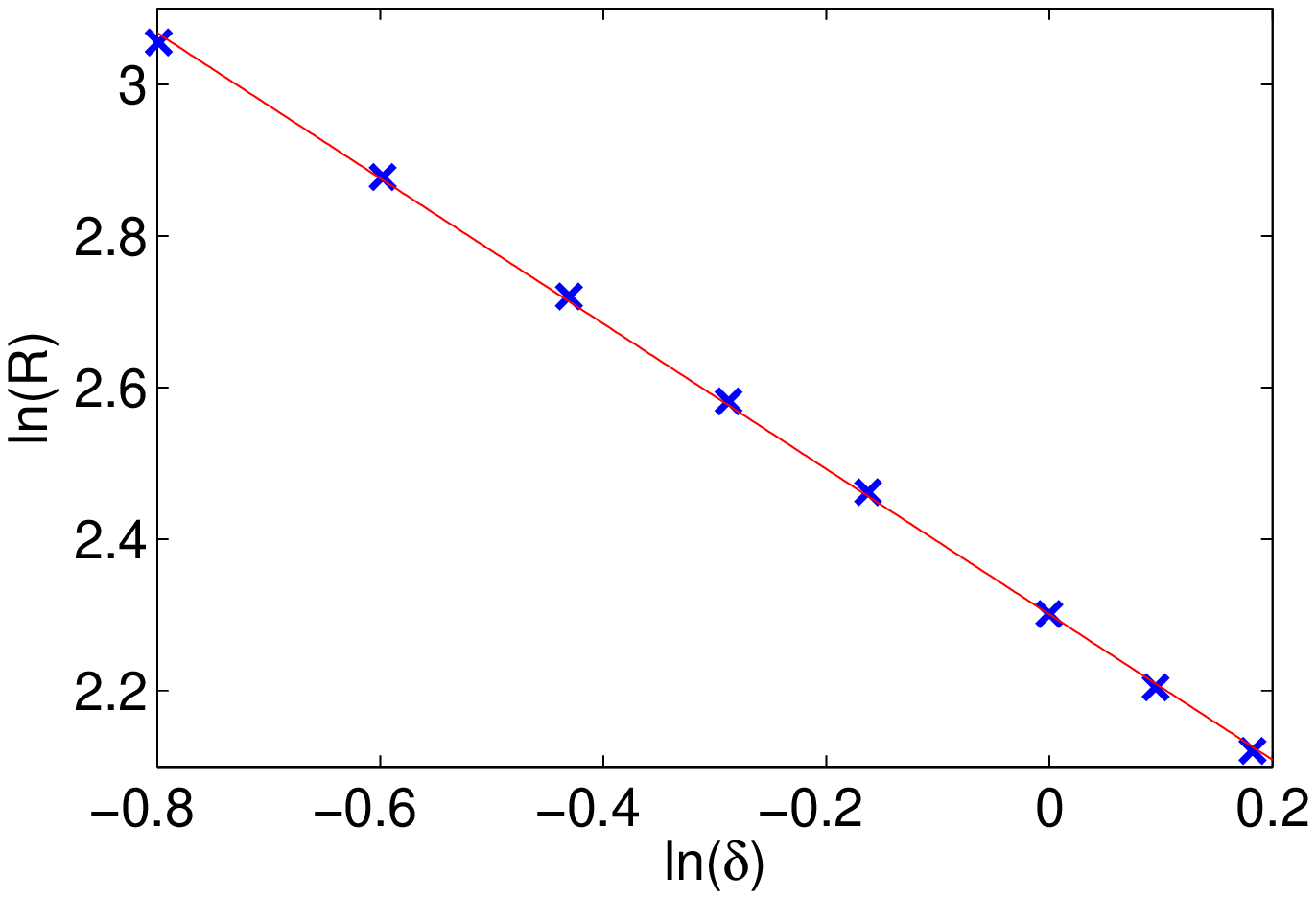}}
\caption{The scaling of the action as a function of the parameters when  $\epsilon=.05$. For the two sweeps, the action is obtained by using the IAMM on Eq. \ref{4dDuff} with $a=1$ and $\delta=1$ while $b$ is swept in (a) and with $a=1$ and $b=20$ while $\delta$ is varied in (b). Since $\epsilon$ is small, we expect the center manifold approximation of Eq. \ref{action_scaling} to be valid. Indeed, the linear best fit in (a) is $\ln(R)=1.9989\ln(b)-2.8417$ while in (b) the equation is $\ln(R)=-.958666\ln(\delta)+2.3008$. The slope of both of these lines are close to the prediction given by equation \ref{action_scaling}. }\label{duff_center_scale}
\end{center}
\end{figure}
\begin{figure}[h!]
\begin{center}
\subfigure[]{\includegraphics[scale=0.35]{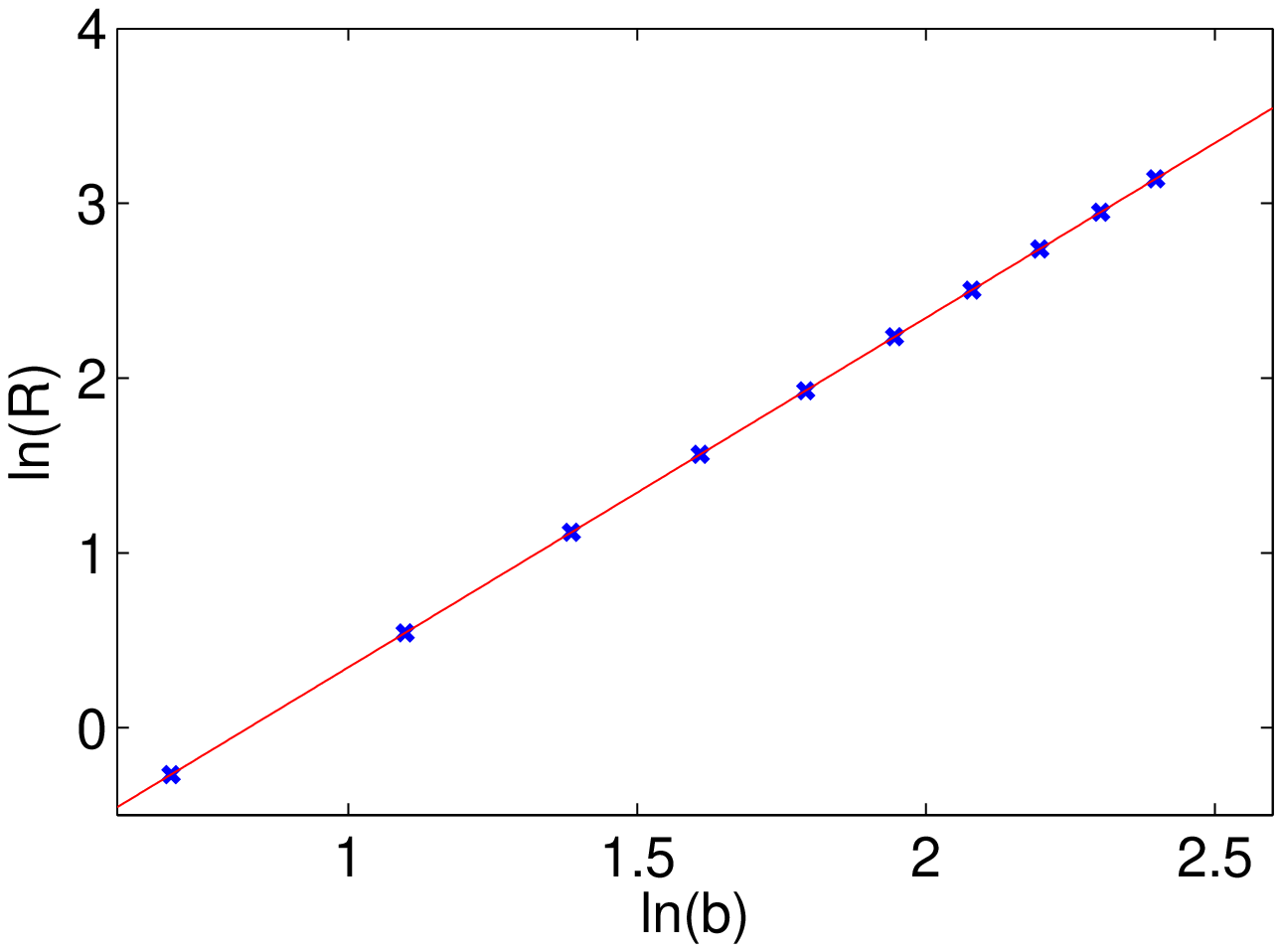}}
\subfigure[]{\includegraphics[scale=0.35]{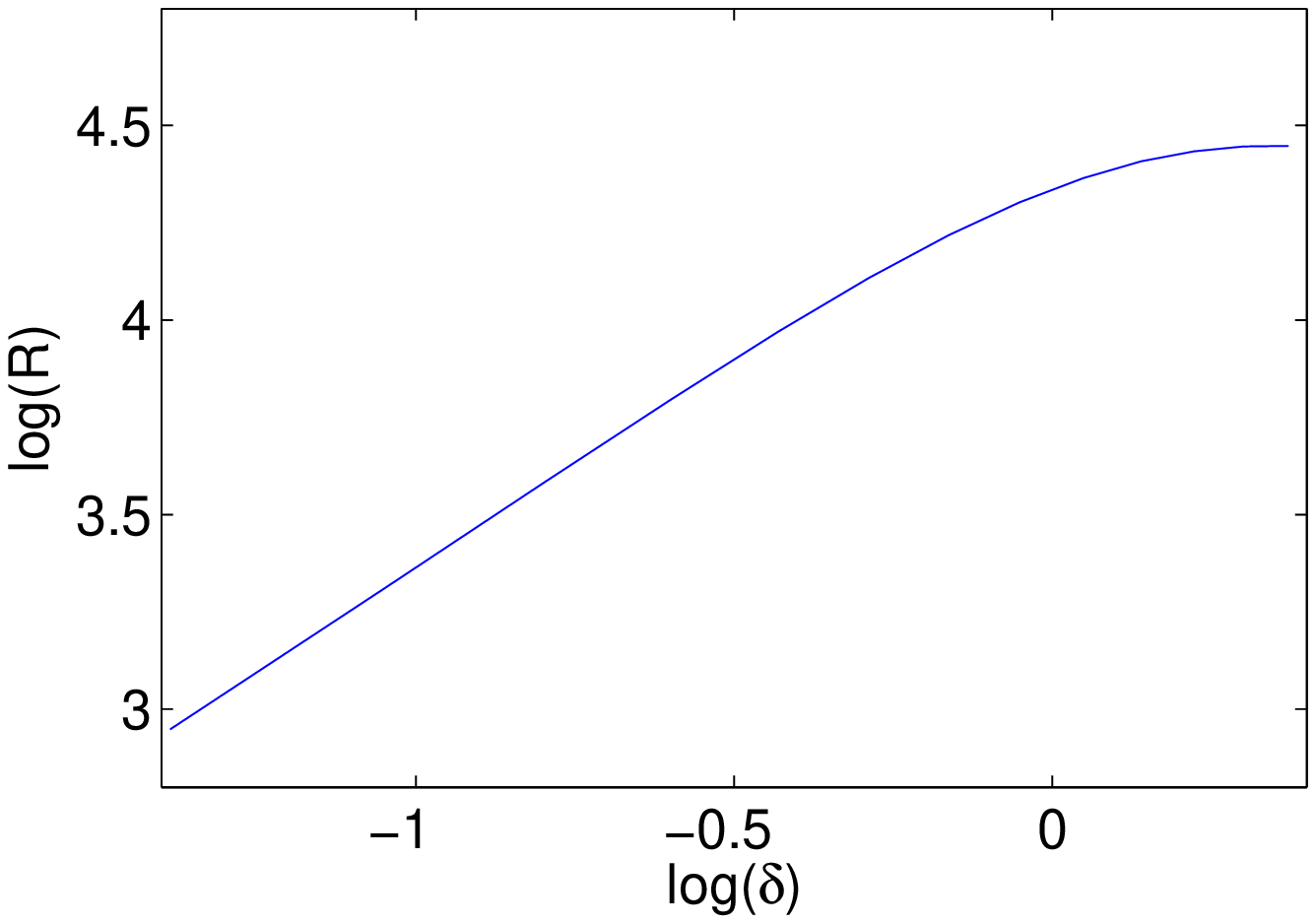}}
\caption{The scaling of the action as a function of the parameters when $\epsilon=1$. For the two sweeps, the nonlinear term $b$ is swept for fixed $\delta=1$ in (a) while $b=20$ is fixed while $\delta$ is varied in (b). The linear best fit in (a) is $\ln(R)=2\ln(b)-1.655$ while in (b) the power law scaling vanishes as $\delta$ is increased.   The second order scaling of the action with respect to $b$ is persistent even away from the center manifold. }\label{duff_scale}
\end{center}
\end{figure}

Since it is not clear if the same relationship will work further away from the
center manifold, we check the scaling with respect to $b$ and $\delta$ when
$\epsilon=1$ in Eq.\ref{4dDuff} in figure \ref{duff_scale}. Interestingly, the leading order scaling of the action with respect to $\beta$ is still $2$, just as it was near the center manifold. Meanwhile, the scaling with respect to $\delta$ is markedly different. Thus, we can predict both near and away from the center manifold how the action will scale as a function of the distance between the two equilibrium points.


\subsection{SIS Epidemic Models}
Next we consider a simple susceptible-infected-susceptible (SIS) epidemic model in the
  general noise case, as considered by \cite{sbdl09}. This model is defined by
  a master equation which describes the probability of fluctuations between a
  susceptible or infected category in a population of $N$ individuals. Assuming $N$ is sufficiently large, we can write a mean field system of equations for the change of the population fractions of susceptible and infected individuals, denoted $x_S$ and $x_I$ respectively,
\begin{equation}
\begin{aligned}
\dot{x}_S & = \mu-\beta x_S x_I + \kappa x_I-\mu x_S \\
\dot{x}_I & = \beta  x_S x_I - \kappa x_I - \mu x_I.
\label{mfsis}
\end{aligned}
\end{equation}
For simplicity, the population in the mean field is assumed constant, i.e. births and deaths are equal. Here $\mu$ is the natural birth and death rate of both the susceptible and infected populations, $\beta$ is the contact rate, and $\kappa$ is the natural recovery rate of the infected population.

\begin{figure}[h!]
\begin{center}
{\includegraphics[scale=0.35]{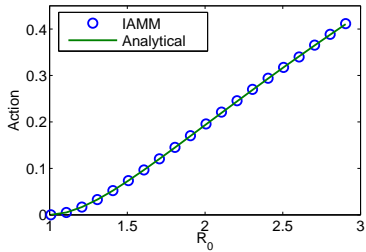}}
\caption{The action as a function of $R_0$ as predicted by an analytical expression and our numerical method. Here, $\kappa = 100$, and $\mu = .2$ while $\beta$ is varied. }\label{1DSISscale}
\end{center}
\end{figure}

Since we have assumed the population fraction of susceptible and infected individuals are conserved, we have $x_I=1-x_S$. Under this constraint, assuming a small random fluctuation of only the infected individuals, we can use a methodology similar to the one described above (and worked out in detail in \cite{sbdl09}) to write a Hamiltonian system,
\begin{equation}H(x_I,p_I) = \beta x_I (1-x_I)(e^{p_I}-1) + (\mu + \kappa)x_I(e^{-p_I}-1).\label{1DSI}\end{equation}
We shall refer to equation \ref{1DSI} as the 1D SIS model equation. The optimal path will extend from the endemic state at $x_I = 1-1/R_0$ to the extinct state $x_I = 0$, where $R_0 =\beta/(\mu+\kappa) $. Using the methods introduced above, the optimal path (at typical parameter values) is given below. Integrating the momentum along this path gives the action, which is proportional to the probability of extinction. In this simple 1D SIS model, the predicted action scaling as a function of $R_0$ is given by solving the Hamiltonian directly for $p_I$ and solving the integral \ref{action:integral} $$R = \int_{1-1/R_0}^0-\text{ln}[R_0(1-x_I)]dx.$$ Since this is one of the rare cases the optimal path can be found analytically, it is a great test case for our method. A comparison of the analytical action to the numerical action is shown in figure \ref{1DSISscale}.

\begin{figure}[h!]
\begin{center}
\subfigure[]{\includegraphics[scale=0.35]{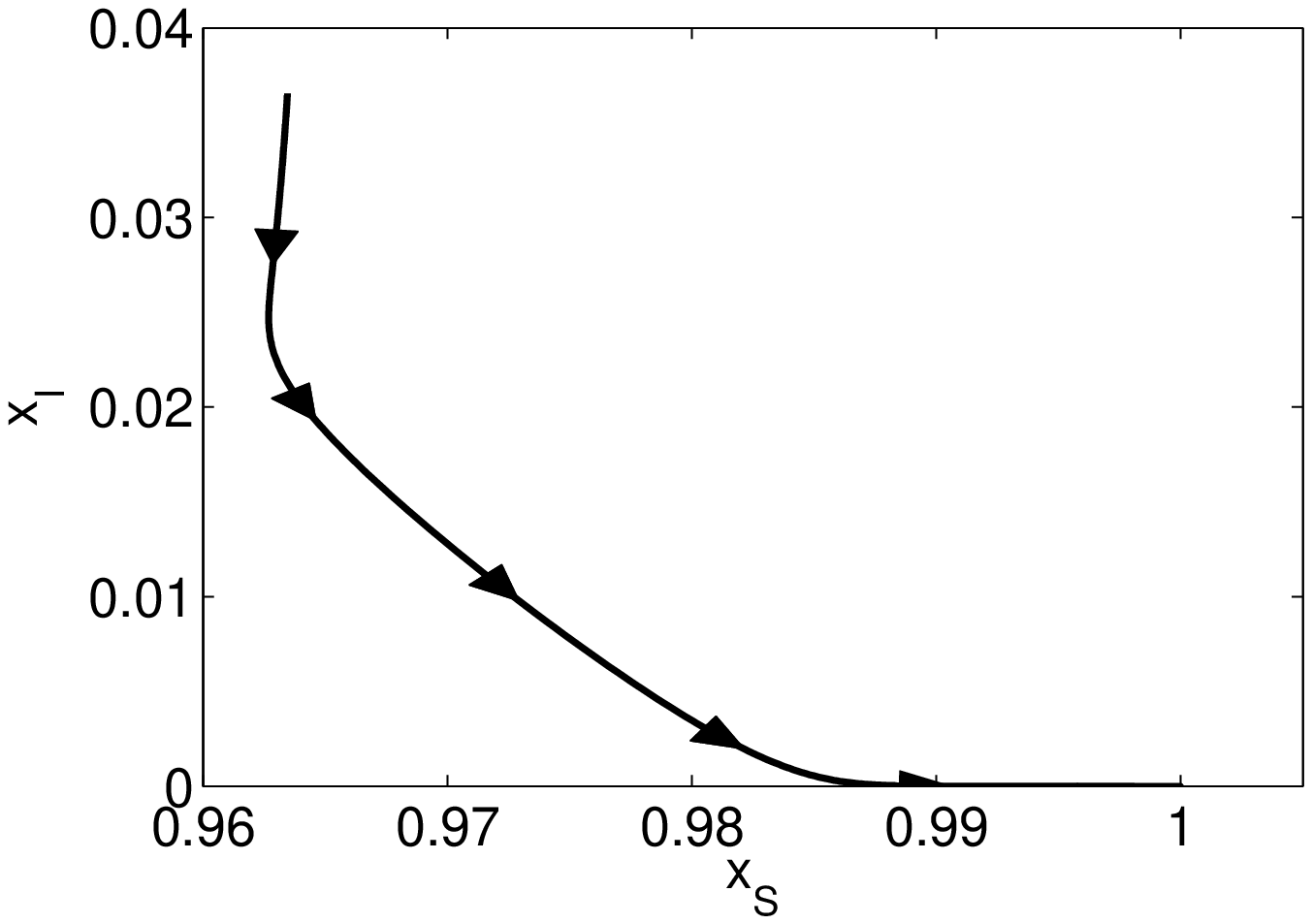}}
\subfigure[]{\includegraphics[scale=0.35]{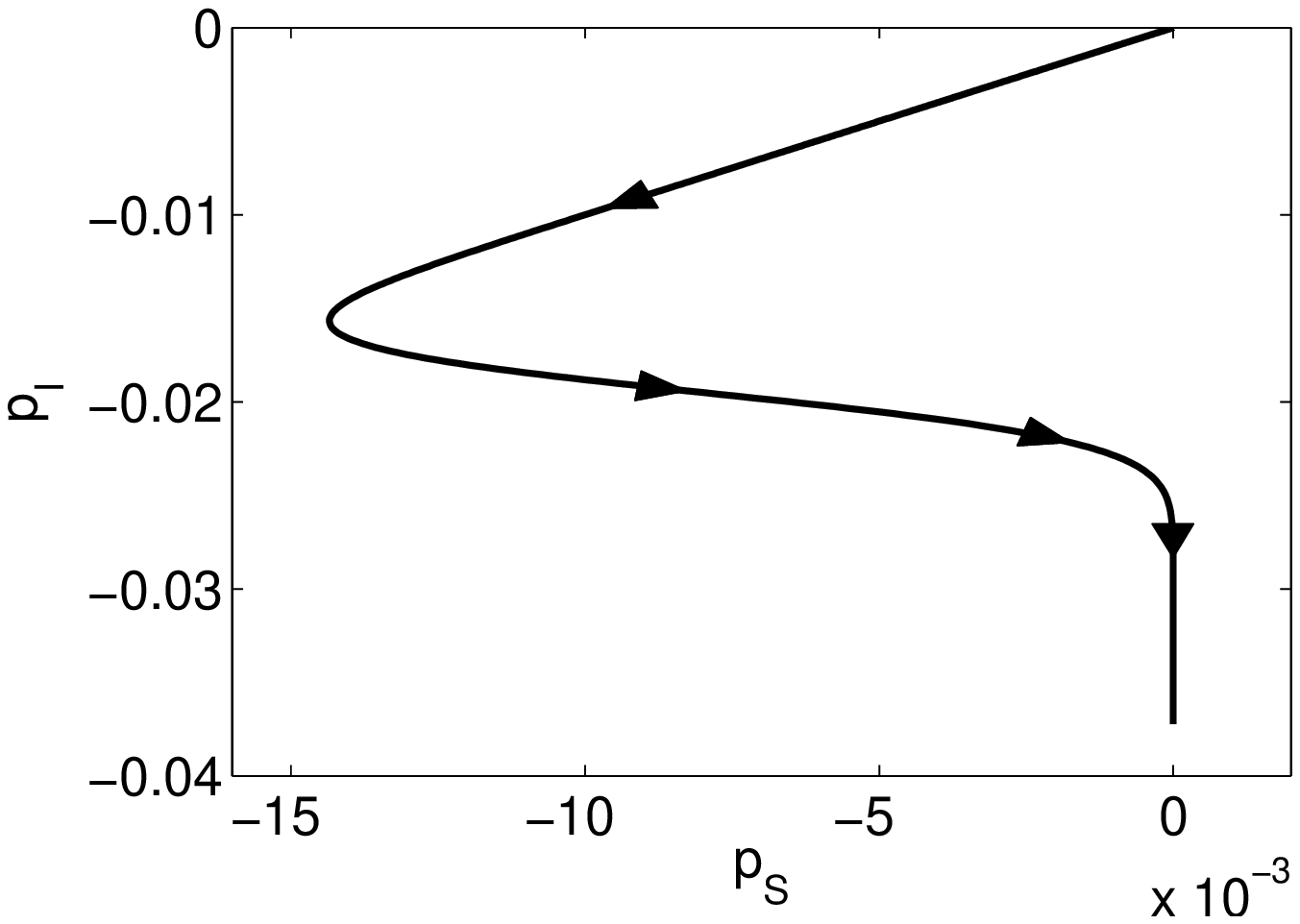}}
\caption{Projection of the optimal path onto (a) population space and (b) momenta space for a sample SIS system with $\beta = 104$, $\kappa = 100$, and $\mu = .2$. }\label{SIS_path}
\end{center}
\end{figure}

In the case of an SIS epidemic model with independent fluctuations on both the susceptible ($x_S$) and infected ($x_I$) populations, the Hamiltonian form of the equations is given by,
\begin{equation}
H(x_S,x_I,p_S,p_I) = \mu (e^{p_S}-1) + \beta x_s x_I (e^{-p_S+p_I}-1) + \kappa x_I (e^{p_S-p_I}-1) + \mu x_S (e^{-p_S}-1) +\mu x_I (e^{-p_I}-1),
\end{equation}
and we refer to this form as the 2D SIS model. The two states of interest for this system are the endemic state, $(1/R_0,1-1/R_0,0,0)$ and the nontrivial extinct state, $(1,0,0,\ln(1/R_0))$. Using the procedure discussed above, we compute the optimal path from the endemic to the extinct state, and show a typical result in figure \ref{SIS_path}. 

\begin{figure}[h!]
\begin{center}
\includegraphics[scale=0.35]{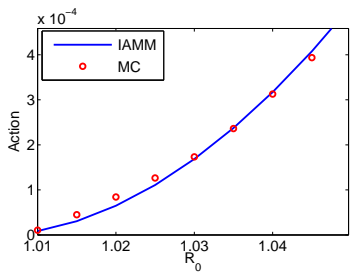}
\caption{A comparison of the action predicted by computing the optimal path to a Monte-Carlo simulation of the original system. Here, $\beta$ is varied while $\kappa=100$ and $\mu=.2$. }\label{SIS_action}
\end{center}
\end{figure}

Instead of comparing the action scaling predicted here to an analytical formulation, we instead compare it to a Monte-Carlo simulation of the master equation for the initial system with a fixed population size of 20,000 individuals. In the Monte-Carlo simulation, a Gillespie algorithm is employed on the SIS model, and from all the simulations, a probability of extinction is computed for several values of $R_0$. From these, a mean extinction time is derived, and the log of this mean extinction time should scale like the action predicted from the action integrals, from Eq. \ref{action:integral}, computed using our approach. Figure \ref{SIS_action} shows a comparison of the two approaches, and good agreement is seen between these two independent methods.

\subsection{Finite Time Lyapunov Exponents}

As demonstrated in \cite{schwartz2011converging}, the optimal path to extinction coincides with ridges, i.e. maximal values, of the finite time Lyapunov exponents (FTLEs). Forgoston and others propose finding FTLE ridges as a method of computing optimal paths \cite{schwartz2011converging, ForgostonBMB, Kessler_2012}. Here, we demonstrate that our approximation to the optimal path does indeed locally maximize the FTLE. 

We proceed by using the methods outlined in
\cite{hall00,hall01,hall02,shlema05} to approximate the FTLEs at points on our
optimal path, and at nearby points transverse to the optimal path. We begin by
picking a point on the optimal path (generated by our method), and on nearby
points some small distance away from the path. 

For the given vector field, we assume we have a flow passing through
  initial point $\bm{x_0}$,
  $\phi:{\cal{R}}^n\rightarrow{\cal{R}}^n$, such that
$\phi_{t_0}^{t_0+T}(\bm{x_0})=\bm{x}(t_0+T;t_0,\bm{x_0})$. The local linear
variation at $\bm(x_0)$ is defined by $\boldsymbol\Delta(\bm{x_0},t_0+T) = \frac{\partial
  \phi_{t_0}^{t_0+T}(\bm{x_0})}{\partial \bm{x_0}}$.

 Using a fourth order Runge-Kutta method, we can integrate all the initial
 points , $\bm{x_0}$, forward in time over a fixed interval, and compute
 the finite time deformation rate of the local coordinates (i.e. Right
 Cauchy-Green Tensor)  $\boldsymbol C(\bm{x_0})
   =\boldsymbol\Delta^T(\bm{x_0},t_0+T)\boldsymbol\Delta(\bm{x_0},t_0+T)$. 
 The maximal eigenvalue $\lambda_{max}$ of $\boldsymbol {C}(\bm{x_0})$  will give the FTLE $\sigma(y,t_i,T)=\frac{1}{T}\ln\sqrt{\lambda_{max}}$.

\begin{figure}[h!]
\begin{center}
\subfigure[]{\includegraphics[scale=0.35]{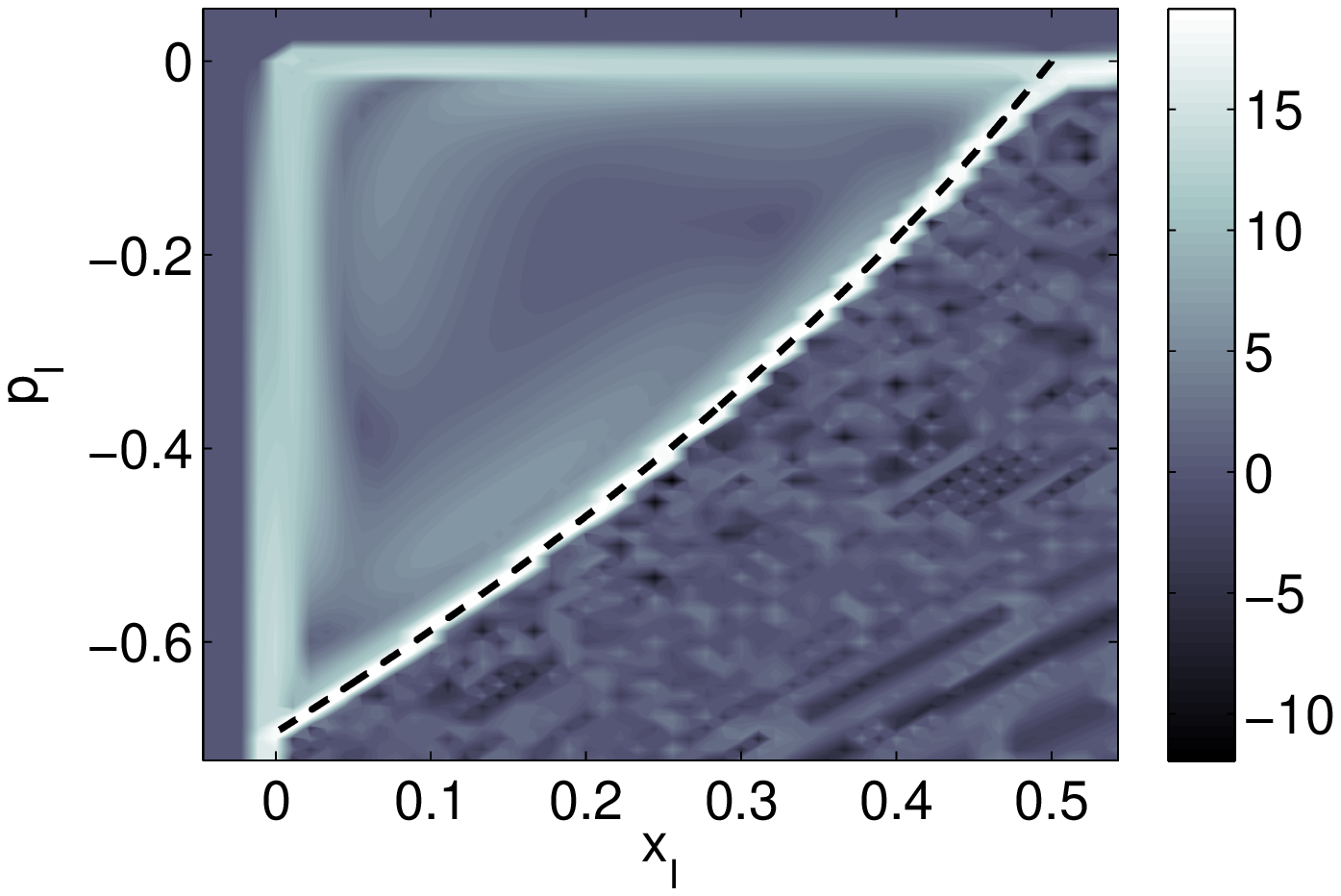}}
\subfigure[]{\includegraphics[scale=0.35]{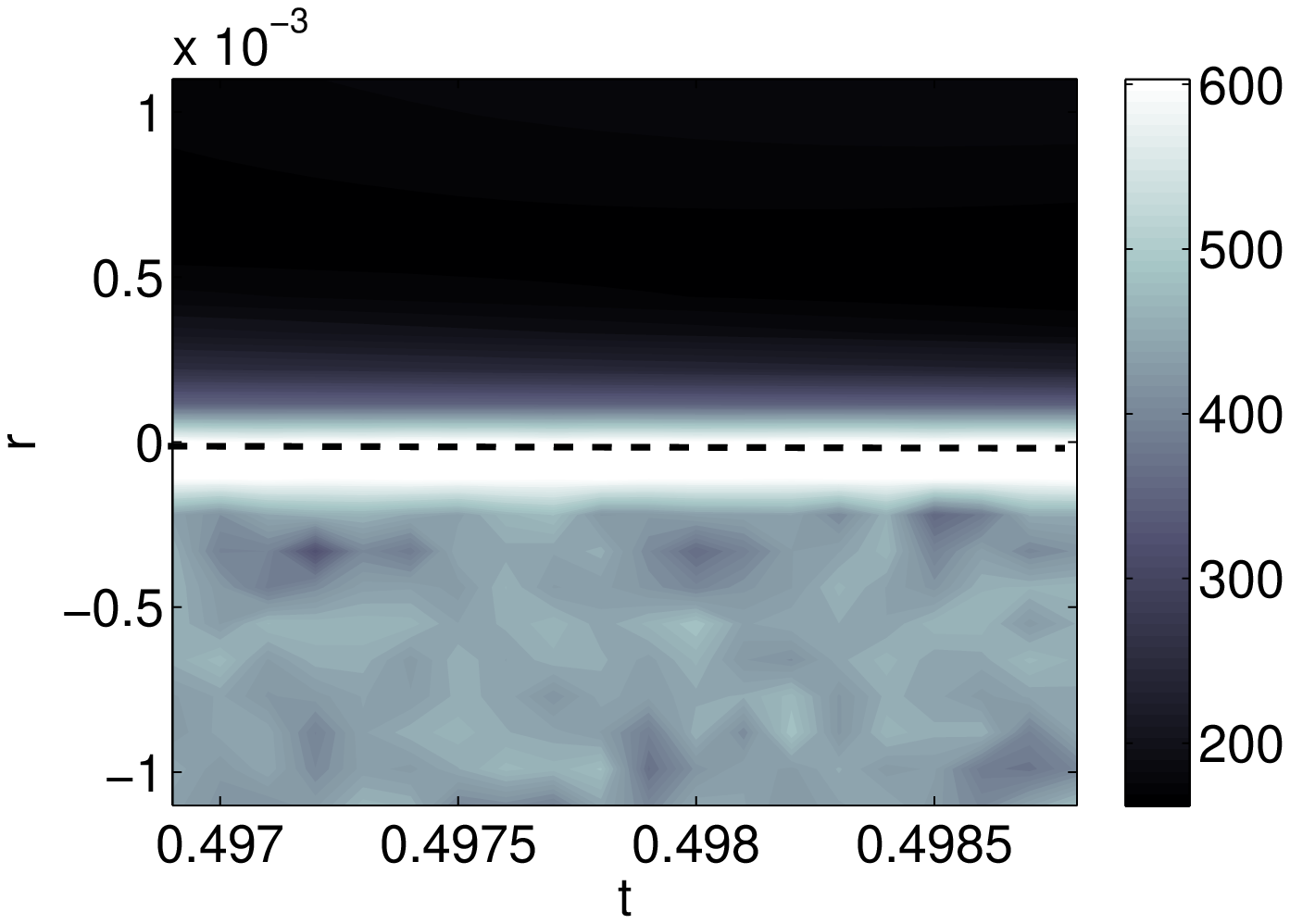}}
\caption{Figure (a) shows the numerically computed optimal path overlaid with the FTLE for the one-dimensional SIS model. The optimal path (dashed line) occurs precisely along a ridge of the FTLE.  Figure (b) shows a single slice of the FTLE computation for the two-dimensional SIS model, where $r$ is the distance from the optimal path in the transverse direction along a unit vector, and $-r$ represented a distance in the antipodal direction. For both paths, the parameter values from Figure \ref{SIS_path} are used.}\label{Lyap1}
\end{center}
\end{figure}

Consider the case of the 1D SIS extinction model from above. Here, we will use a path computed above and compare the values near this path to the local FTLEs. In Figure \ref{Lyap1}a we plot the FTLEs over a square domain, and show that a local maximum (ridge) is attained precisely where the computed optimal path predicts.  

In higher dimensions, the maximal Lyapunov exponent is still exhibited along
the optimal path~\cite{schwartz2011converging}. For the 2D SIS model, note that the optimal path exists in
four dimensions. Thus, the transverse direction is the set of points obtained by rotating a
normal vector to the path around two Euler angles (which forms a 3D sphere for a given radius $r$). To illustrate FTLEs in this higher dimension framework, we must consider the
``shell'' around the initial starting point (and orthogonal to the path) at a
fixed radius, and then compute FTLEs on this shell all along the optimal path. We expect that the maximal FTLE will occur along the optimal path, relative to nearby points on the transverse sphere. To illustrate this, we define $u_{\perp}(t)$ as a unit vector orthogonal to the tangent vector of the optimal path at a given time, i.e.  $u_\perp(t) = (\dot{x}_I(t),-\dot{x}_S(t),0,0)/|(\dot{x}_I(t),-\dot{x}_S(t),0,0)|$, and then examine the FTLE as a function of the two Euler angles for a given $r$. Figure \ref{Lyap1}b shows just one cross section (over a short time interval), obtained by setting both rotation angles to either $0$ or the anitipodal angle $\pi$, of this high dimensional object, as a function of $r$, and demonstrates that the maximal FTLE is, indeed, along the optimal path.

\section{ Time Delayed SDEs}

One advantage of the IAMM is that it allows the solution of stochastic delay-differential equations of the form,
\begin{equation}
\dot{\vx}(t) = \boldsymbol{f}(\vx(t),\vx(t-\tau)) + \boldsymbol{G}(\vx(t))\boldsymbol{\xi}(t).
\label{general:dde}
\end{equation}
Schwartz et. al \cite{Schwartz2012} have demonstrated that the methodology introduced in section 2 can be adapted to write this system as a Hamiltonian system,
\begin{equation}
H(\bm{x},\bm{x}_{\tau},\vl) =
\frac{(\boldsymbol{G}^2(\boldsymbol{x})\boldsymbol{p})\cdot\boldsymbol{p}}{2}+\boldsymbol{p}\cdot
\boldsymbol{f}(\vx,\vx_\tau),\label{hamiltonian:delay}
\end{equation}
where $\vx_\tau=\vx(t-\tau)$. The equations of motion are given by,
\begin{equation}
\begin{aligned}
\dot{\vx} &= \frac{\partial H}{\partial \vp}(\bm{x},\vx_\tau,\bm{p}) \\ \dot{\vp} &= -\frac{\partial H}{\partial \vx}(\bm{x},\vx_\tau,\bm{p}) - \frac{\partial H}{\partial \vx_\tau}(\bm{x}(t+\tau),\vx(t),\bm{p}(t+\tau)).
\label{delay:eom}
\end{aligned}
\end{equation}
Note the appearance of both delay and advance terms in \ref{delay:eom}. Because of the appearance of the delay term, the Hamiltonian is no longer time invariant, and unlike in the previous examples, where $H(\vx,\vp)=0$, the zero-energy condition is not conserved.

We shall consider a one dimensional test case where $f(x,x_\tau) = x(1-x)-\gamma x_\tau$, where the steady states are given by $x_A = 1-\gamma$ and $x_B = 0$, for $0\leq\gamma\leq 1$. Again, we will assume additive noise $\boldsymbol{G}=1$, and derive the Hamiltonian system,
\begin{equation}
H(x,x_\tau,p) = \frac{p^2}{2}+p(x(1-x)-\gamma x_\tau),
\end{equation}
and the corresponding equations of motion,
\begin{equation}
\begin{aligned}
\dot{x} &= x(1-x)-\gamma x_\tau +p\\
\dot{p} &= -p(1-2x) + \gamma p(t+\tau).
\label{delay:eom1}
\end{aligned}
\end{equation}
The IAMM needs only a few minor adjustments to compute these paths numerically. Primarily, the presence of the delay and advance terms will add additional entries into our linear system Eq. \ref{discrete:equations}. Our method uses a non-uniform timestep, and so $x(t_k-\tau)$ may not coincide exactly with one of our points $x_k$ at time $t_k$. To overcome this, we can use Lagrange interpolation on the closest four points $x_{j-2},x_{j-1},x_j,x_{j+1}$ such that $t_{j-2}<t_{j-1}<t_k-\tau<t_j<t_{j+1}$. Since we keep the time domain fixed the $j$ needed for each $t_k$ can be easily computed before the iterative scheme is started, and fewer or more terms can be used in the Lagrange interpolation scheme depending upon desired accuracy. Further, if $t_k-\tau<T_\epsilon$ or  $t_k+\tau>T_\epsilon$, i.e. the delay or advance terms fall outside of our numerical domain, then we can set $x_k = x_A$ or $x_k = x_B$ respectively.

\begin{figure}[h!]
\begin{center}
\includegraphics[scale=0.35]{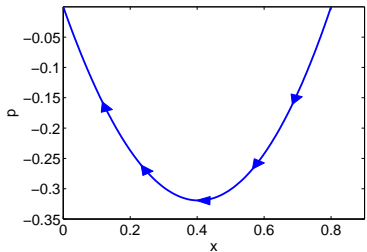}
\caption{The optimal path for eq. \ref{delay:eom1} with $\gamma=.2$ and $\tau=.5$.}\label{delaypath}
\end{center}
\end{figure}

\begin{figure}[h!]
\begin{center}
\includegraphics[scale=0.35]{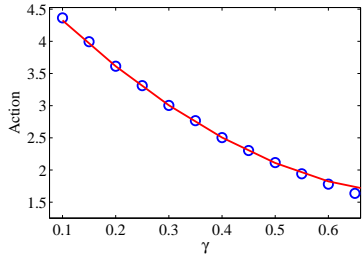}
\caption{Scaling of the action along the optimal path as a function of $\gamma$. Here the solid line indicates the IAMM predication, and the circles represent a Monte-Carlo simulation of 1000 runs. }\label{delayscaling}
\end{center}
\end{figure}

To demonstrate the effectiveness of the IAMM in solving these stochastic delay problems, we show sample paths in figure \ref{delaypath} and compare the scaling of the action versus a Monte-Carlo simulation of the stochastic delay difference equation in figure \ref{delayscaling}, and note the good agreement.

\section{Discussion}

We have considered the  problem of finding the trajectory in stochastic
  dynamical systems that optimizes the
  probability of switching between two states, or causes one or
  more components to go extinct. In computing such a trajectory, called the
  optimal path, we needed to consider a numerical technique which could  solve
  a Hamiltonian system with asymptotic boundary conditions in time. We
 have developed a numerical method, which we call the iterative action minimizing method (IAMM), for finding the
optimal path of transition between two steady states in stochastic dynamical
systems. This method is ideal for systems which can be written as two-point
boundary value  problems governed by Hamiltonian systems. We have validated the IAMM by presenting
a variety of problems of interest, and have compared the numerical results
with either analytic results or Monte-Carlo simulations of full stochastic
systems.

As demonstrated here, the IAMM method is robust enough to be applicable to
a variety of different types of problems, including continuous SDE systems,
such as the Duffing equation, discrete epidemic models of finite population
size, such as the SIS model,
and stochastic delay differential equations, in which the deterministic
problem is infinite dimensional. The methodology is
straightforward enough to generalize to higher dimensions, in contrast to
other commonly used methods, such as the shooting method, which is a major
advantage of the IAMM. 

The primary limitations of this method are scaling issues in very high
state space dimensions, and the finesse required in picking an initial
guess that guarantees convergence, both of which are typical of iterative
methods of quasi-Newton type. In the limit of small noise or large system
size, however,  due to the robustness and ease of
generalization to complex and high dimensional dynamical systems, the
  method offers a considerable advantage over simulating large systems, or systems
  which require many Monte Carlo runs to generate statistics of the
  transitions paths. As a result, we expect 
 this method
will be useful in efficiently solving a large variety of optimal transition problems in
the field of stochastic dynamical systems. 


\section{Acknowledgments}
The authors gratefully acknowledge the Office of Naval Research for their
support under  N0001412WX20083, and support of the NRL Base Research Program
N0001412WX30002. Brandon Lindley is currently an NRC Postdoctoral Fellow. We thank Lora Billings for providing the Monte Carlo data used in figures \ref{SIS_action} and \ref{delayscaling}, and Eric Forgoston for a preliminary reading of this manuscript.



\end{document}